\documentclass[12pt,reqno]{amsart}
\usepackage{amsmath,amssymb,amsfonts,amscd,latexsym,amsthm,mathrsfs,verbatim}
\usepackage[unicode]{hyperref}
\usepackage{hypbmsec}
\DeclareMathOperator*{\CT}{CT}
\renewcommand{\geq}{\geqslant}
\renewcommand{\leq}{\geqslant}

\begin{document}

\title{$q$-rious and $q$-riouser}

\author{S.~Ole Warnaar}

\author{Wadim Zudilin}

\dedicatory{To our friend Dick Askey}
\date{September 2019}

\maketitle

Dick Askey is known not just for his beautiful mathematics and his
many amazing theorems, but also for posing numerous interesting and 
important open problems.
Dick being Dick, these problems are hardly ever isolated, and often 
intended to demonstrate the unity of analysis, number theory and
combinatorics.
On this ocassion we wish to take the reader down the rabbit hole created by
one such problem, published as \emph{Advanced Problem 6514} by the
American Mathematical Monthly in April 1986 \cite{Askey86}.
Dick's inspiration for the problem was derived from the Macdonald--Morris
constant term conjecture for the root system $\mathrm{G}_2$ 
\cite{Macdonald82,Morris82}
as well as much earlier work of P.~Chebyshev \cite{Tchebichef82} 
and E.~Landau \cite{Landau85} on the integrality of factorial ratios.
Problem 6514 asks for a proof of the integrality of
\[
A(m,n)=
\frac{(3m+3n)!\,(3n)!\,(2m)!\,(2n)!}{(2m+3n)!\,(m+2n)!\,(m+n)!\,m!\,n!\,n!}
\]
for all non-negative integers $m$ and $n$.

There are multiple reasons --- some of them very deep, see e.g.,
\cite{Bober09,Rodriguez-Villegas05,Soundararajan19a,Soundararajan19b}
--- for wanting to classify integer-valued factorial ratios such as
Chebyshev's
\[
C(n)=\frac{(30n)!\,n!}{(15n)!\,(10n)!\,(6n)!}\,.
\]
Given a particular such ratio, integrality can always be verified
by computing the $p$-adic order of the factorials entering the quotient.
This is exactly what all eight solvers of Problem 6514 did.
Such a verification, however, provides little insight into which ratios 
are integral and which ones are not, and from the editorial comments
to the problem it is clear that Dick would have liked to see other types
of solutions too.
Indeed, it is remarked that
\begin{quote}
[\;\;] the editor [read: Dick Askey] feels there is still room 
for other methods, involving perhaps combinatorial interpretations
or manipulation of generating functions. In this particular case, 
the proposer remarks that $A(m,n)$ should be the constant term of 
the Laurent polynomial 
\begin{multline*}
\quad\qquad \big((1-x)(1-1/x)(1-y)(1-1/y)(1-y/x)(1-x/y)\big)^m \\[1mm]
\times \big((1-xy)(1-1/xy)(1-y/x^2)(1-x^2/y)(1-y^2/x)(1-x/y^2)\big)^n.\quad
\end{multline*}
\end{quote}
Incidentally, L.~Habsieger \cite{Habsieger86} and D.~Zeilberger 
\cite{Zeilberger87} both proved the $\mathrm{G}_2$ Macdonald--Morris
constant term conjecture shortly after Dick Askey posed his problem.
The submission dates of their respective papers (the 12th of May and the
2nd of June 1986) were well within the deadline of the 31st of August
for submitting solutions to Problem 6514 to the Monthly.
In fact, in the acknowledgement of his paper Zeilberger thanks
Dick Askey for ``rekindling his interest in the Macdonald conjecture'',
so maybe he should belatedly be considered the 9th solver of Askey's problem.

The height of a factorial ratio is the number of factorials in the
denominator minus the number of factorials in the numerator, so that the
height of $A(m,n)$ is two whereas the height of $C(n)$ is one. 
A one-parameter family of height-$k$ factorial ratios 
\[
F(n)=\frac{(a_1\,n)!\cdots (a_{\ell}\, n)!}
{(b_1\, n)!\cdots (b_{k+\ell}\, n)!}
\]
is balanced if $a_1+\dots+a_{\ell}=b_1+\dots+b_{k+\ell}$.
All balanced, integral, height-one factorial ratios $F(n)$
were classified in 2009 by J.~Bober~\cite{Bober09}.
In relation to this classification we should mention 
F.~Rodriguez-Villegas' observation \cite{Rodriguez-Villegas05} 
that if $F(n)$ is a balanced, height-one factorial ratio
then the hypergeometric function $\sum_{n\geq 0} F(n) z^n$ is algebraic 
if and only if $F(n)$ is integral. This observation was 
key to Bober's proof, allowing him to use the Beukers--Heckman
classification \cite{BH89} of $_nF_{n-1}$ hypergeometric functions with
finite monodromy group. A proof not reliant on the
Beukers--Heckman theory was recently found by 
K.~Soundararajan \cite{Soundararajan19a}.
By extending his method he also obtained a partial
classification in the height-two case \cite{Soundararajan19b}.

Despite the availability of the number-theoretic, $p$-adic approach to 
factorial ratios, the question of integrality is very interesting from
a purely combinatorial point of view.
The simplest example is of course provided by the height-one binomial
coefficients
\[
\frac{(m+n)!}{m!\,n!},
\]
whose integrality can be established combinatorially 
(as well as probabilistically, algebraically, etc.) with little effort.
However, to the best of our knowledge, no combinatorial proof is
known of the integrality of Chebyshev's $C(n)$.

A related open problem arises from our joint work \cite{WZ11} from 2011.
In \cite{WZ11} we observed that if each factorial $m!$ in an integral 
factorial ratio is replaced by a $q$-factorial
\[
[m]!=[m]_q!=\prod_{i=1}^m\frac{1-q^i}{1-q},
\]
then the resulting $q$-factorial ratio is a polynomial 
with non-negative integer coefficients.
The polynomiality and integrality parts are trivial but the positivity ---
which was referred to in \cite{WZ11} as `$q$-rious positivity' ---
is completely open.
The only (irreducible) cases that are proven are the three 
two-parameter families of height one given by
\[
\frac{[m+n]!}{[m]!\,[n]!},\qquad
\frac{[2m]!\,[2n]!}{[m]!\,[n]!\,[m+n]},\qquad
\frac{[m]!\,[2n]!}{[2m]!\,[n]!\,[n-m]!}\quad (m\leq n),
\]
where the first family corresponds to the $q$-binomial coefficients
and the second family to the $q$-super Catalan numbers.
In the $q$-case no arithmetic approach is available, and given the lack of 
combinatorial methods to deal with integrality, a combinatorial
approach to $q$-rious positivity seems hopeless.\footnote{There are of course
countless methods to show that the $q$-binomial coefficients have non-negative
integer coefficients, but no combinatorial interpretation of the $q$-super 
Catalan numbers is known. In fact, not even a combinatorial 
interpretation of the ordinary super Catalan numbers is known.}
Perhaps the most tractable problem is to analytically prove, along the 
lines of \cite{WZ11}, the positivity of the known two-parameter families of 
height two, such as
\[
A_q(m,n)=
\frac{[3m+3n]!\,[3n]!\,[2m]!\,[2n]!}
{[2m+3n]!\,[m+2n]!\,[m+n]!\,[m]!\,[n]!\,[n]!}\in\mathbb Z[q]
\]
and
\[
C_q(m,n)=\frac{[6m+30n]!\,[n]!}
{[3m+15n]!\,[2m+10n]!\,[m]!\,[6n]!}\in\mathbb Z[q].
\]
For the first family, which is the $q$-analogue of $A(m,n)$,
it is known that \cite{Cherednik95,Habsieger86,Zeilberger87}
\begin{multline*}
A_q(m,n) \\=\CT\limits_{x,y}\Big[
\big(x,q/x,y,q/y,y/x,qx/y;q\big)_m
\big(xy,q/xy,y/x^2,qx^2/y,y^2/x,qx/y^2;q\big)_n\Big],
\end{multline*}
where $(a_1,\dots,a_k;q)_n:=\prod_{i=1}^k \prod_{j=1}^n (1-a_i q^{j-1})$.
This interpretation as a $\mathrm{G}_2$ constant term gives little 
insight into the positivity of the coefficients.
It would appear that the second two-parameter family has not
occurred before. For $q=1$ it arose earlier this year in the (partial)
classification of height-two factorial ratios by 
Soundararajan~\cite{Soundararajan19b} mentioned above.
It should be noted that if one were to prove the $q$-rious positivity of
$C_q(m,n)$ then this immediately would imply the positivity of the
$q$-analogue of Chebyshev's factorial ratio since
\[
C_q(n)=\frac{[30n]!\,[n]!}{[15n]!\,[10n]!\,[6n]!}=C_q(0,n).
\]

\end{document}